\documentclass[11pt, makeidx]{amsart}

\usepackage{enumerate}
\usepackage{amssymb,amsfonts,amsthm}
\usepackage{srcltx}

\usepackage{float}

\usepackage[OT2,OT1]{fontenc}
\newcommand\cyr{%
\renewcommand\rmdefault{wncyr}%
\renewcommand\sfdefault{wncyss}%
\renewcommand\encodingdefault{OT2}%
\normalfont
\selectfont}
\DeclareTextFontCommand{\textcyr}{\cyr}

\usepackage{color}
\makeindex

\makeatletter
\def\sssub{\@startsection{paragraph}{4}}

\renewcommand\paragraph{\@startsection{paragraph}{4}{\z@}{1.25ex}{0.0001pt}{\normalfont\normalsize\em}}

\numberwithin{paragraph}{subsubsection}

\newcommand{\R}{{\mathbb R}}

\newcommand{\N}{\mathbb{N}}

\newcommand{\Con}{{\mathrm{Con}}}

\newcommand{\dist}{{\mathrm{dist}}}

\newtheorem{theorem}{Theorem}

\theoremstyle{definition}
\newtheorem{df}[theorem]{Definition}

\newcommand{\iv}{^{-1}}

\newenvironment{dedication}
    {\vspace{3ex}\begin{quotation}\begin{center}\begin{em}}
    {\par\end{em}\end{center}\end{quotation}}
\begin{document}

\title[On groups with locally compact asymptotic cones]{On groups with locally compact asymptotic cones} \author{Mark Sapir}\thanks{The research was supported in part by NSF and BSF
grants} \address{Department of Mathematics, Vanderbilt
University, Nashville, TN 37240, U.S.A.} \email{m.sapir@vanderbilt.edu}
\subjclass[2000]{{Primary 20F65; Secondary 20F69, 20F38, 22F50}}
\keywords{asymptotic cone, locally compact, group}

\maketitle
\begin{dedication}
Dedicated to Stuart Margolis' 60s birthday
\end{dedication}
\bigskip

\begin{abstract} We show how a recent result of Hrushovsky \cite{H} implies that if an asymptotic cone of a finitely generated group is locally compact, then the group is virtually nilpotent.
\end{abstract}

\bigskip

Let $G$ be a group generated by a finite set $X$.
Then $G$ can be considered as a metric space where $\dist(g,h)$ is the length of a shortest word on $X\cup X\iv$ representing
$g\iv h$. Let $\omega$ be a non-principal  ultrafilter on $\N$, i.e. a function from the set of subsets $P(\N)$ to $\{0,1\}$ with $\omega(A\cup B)=\omega(A)+\omega(B)$ if $A\cap B=\emptyset$, $\omega(\N)=1$ and $\omega(A)=0$ for every finite set $A$.

Note that if $A\subseteq \N$, then $\omega(A)+\omega(\N\setminus A)=1$ because $\N$ is a disjoint union of $A$ and $\N\setminus A$. For every $A\subseteq B\subseteq \N$, we have  $\omega(A)
\le \omega(B)$ because $B$ is a disjoint union of $B\setminus A$ and $A$. For every two subsets $A,B\subseteq \N$, $\omega(A\cup B)\le \omega(A)+\omega(B)$ because $A\cup B$ is a disjoint union of $A\setminus B\subseteq A$ and $B$. Also if $\omega(A)=\omega(B)=1$, then $\omega(A\cap B)=1$ because $A\cap B=\N\setminus ((\N\setminus A)\cup (\N\setminus B))$.  For every property $P$ of natural numbers, we say that $P$ is true $\omega$-almost surely if the set $S$ of numbers with property  $P$ satisfies $\omega(S)=1$. Thus is $\omega$-almost surely every natural number satisfies $P$ and $\omega$-almost surely every natural number satisfies $Q$, then $\omega$-almost surely every natural number satisfies $P$ and $Q$. We are going to use this property of ultrafilters several times later without reference.

For every sequence of real numbers $r_i\ge 0$ one can define the limit $\lim_\omega r_i$ as the (unique) number $r\in \R\cup\{\infty\}$ such that for every $\epsilon$ we have $|r_i-r|<\epsilon$ $\omega$-almost surely for all $i\in \N$ (see \cite{D}).

Choose a sequence of \emph{scaling constants} $d_n>0$ such that $\lim_\omega d_n =\infty$ and define a pseudo-metric $\overline\dist$ on the Cartesian power $G^\N$ as $$\overline\dist((u_i),(v_i))=\lim_\omega \frac{\dist(u_i,v_i)}{d_i}.$$ Let $G^\N_b$ be the connected component of the sequence $\bar 1=(1,1,\ldots)\in G^\N$. Thus $G^\N_b$ consists of all the elements from $G^\N$ at finite distance from $\bar 1$. Let $\sim$ be the equivalence relation $(u_i)\sim (v_i)$ if and only if $\overline\dist((u_i),(v_i))=0$. Then $\overline\dist$ induces a metric $\dist^\omega$ on $G^\N_b/\sim$, and $G^\N_b/\sim$ with this metric is called the \emph{asymptotic cone} of $G$ corresponding to the ultrafilter $\omega$ and the sequence of scaling constants $(d_n)$, denoted $\Con^\omega(G,(d_i))$ (for more details see \cite{D,DS,SBMS}). Elements of $\Con^\omega(G,(d_n))$ corresponding to sequences $(y_i)\in G^\N$ will be denoted by $(y_i)^\omega$. Note that $G^\N_b$ is a group which acts on $\Con^\omega(G,(d_i))$ transitively by isometries: $(g_i)\cdot (u_i)^\omega=(g_iu_i)^\omega$ (see \cite{D,SBMS}).

Asymptotic cones were first explicitly used by van den Dries and Wilkie \cite{DW} to prove the celebrated theorem of Gromov \cite{G1} that every finitely generated group of polynomial growth is virtually nilpotent (for more applications see \cite{G2,D,DS,SBMS}). One can easily deduce from \cite{DW} that if all asymptotic cones of a finitely generated group are proper (i.e. all closed balls are  compact), then $G$ has polynomial growth and hence is virtually nilpotent.

A finitely generated group can have several non-homeomorphic asymptotic cones \cite{TV}. In fact there are finitely generated groups with continuum non-$\pi_1$-equivalent asymptotic cones \cite{DS}. On the other hand, every virtually nilpotent finitely generated group has only one asymptotic cone up to isometry by Pansu \cite{Pan}.
In \cite{DW}, van den Dries and Wilkie asked whether local compactness of just one asymptotic cone of $G$ implies that $G$ is virtually nilpotent.

F. Point proved \cite{P} that if one asymptotic cone of $G$ is locally compact and has finite Minkovsky dimension, then $G$ is virtually nilpotent.

Here we will show how to use a recent result of Hrushovsky \cite[Theorem 7.1]{H} to answer the question of van den Dries and Wilkie affirmatively. To formulate Hrushovsky's result we need the notion of an approximable subgroup.

\begin{df}\label{d:a} Let $k\in \N$. A finite subset
$X$ of a group $G$  is said to be a $k$-approximate subgroup if $1 \in X$, $X = X\iv,$
and $XX$ is contained in the union of $k$ left
cosets $gX$ of $X$. (Hrushovsky uses right cosets but it does not affect the results because $x\mapsto x\iv$ is an anti-isomorphism for every group.)
\end{df}

\begin{theorem}[Hrushovsky \cite{H}] \label{t:71}
Let
$G$ be
a finitely generated group, $k\in N$.
Assume that $G$ has a collection $\{F_i, i\in \N\}$ of
$k$-approximate subgroups such that any finite subset $F \subset G$ is contained in a one of the $F_i$. Then $G$ is virtually nilpotent.
\end{theorem}

\begin{theorem} \label{t:1} If one asymptotic cone of a finitely generated group $G$ contains a compact closed ball of positive radius, then $G$ is virtually nilpotent.
\end{theorem}

\proof Suppose that a closed ball of radius $\epsilon>0$ in $C=\Con^\omega(G, (d_i))$ is compact. Since $C$ admits a transitive group of isometries, we can assume that this is the ball  $B_C(\epsilon,\bar 1)$ of radius $\epsilon$ around $\bar 1$. Since $B_C(\epsilon,\bar 1)$ is compact, it is contained in the union of a finite number of balls $B_C(\frac{\epsilon}{4}, x_i)$, $x_i\in B_C(\frac78\epsilon,\bar 1)$, $i=1,\ldots,k$.  Let $x_i=(g_{i,j})^\omega$, $i=1,\ldots,k$, $j\in \N$. Note that by the definition of asymptotic cone, for every $i=1,\ldots,k$, $\dist(g_{i,j},1)< \epsilon d_j$ $\omega$-almost surely (since $\frac78<1$).

For every $j\in \N$ let $2p_j$ be the smallest even natural number which is greater than $\epsilon d_j$.
Then $|2p_j-\epsilon d_j|< 2$ for every $j$, whence $\lim_\omega p_j=\infty.$

For every $i=1,\ldots,k$ and $j\in \N$ consider the balls $B_G(2p_j,1)$ and $B_G(p_j,g_{i,j})=g_{i,j}B_G(p_j,1)$ in $G$ (viewed as a metric space).

We claim that $\omega$-almost surely for every $j$, the ball $B_G(2p_j,1)$ is covered by the balls $B_G(p_j,g_{i,j})$, $i=1,\ldots,k$. Indeed, let $S$ be the set of natural numbers $j$ for which this is true. Suppose that  $\omega(S)=0$.  Then $\omega(\N\setminus S)=1$. For each $j\in \N\setminus S$ let $y_j$ be a point in $B_G(2p_j,1)$ which does not belong to any of the balls $B_G(p_j,g_{i,j})$. Consider any sequence $(z_i)\in G^\N$ such that $z_j=y_j$ for all $j\in \N\setminus S$. Since $\omega(\N\setminus S)=1$, $z=(z_i)^\omega$ is at distance at least $\lim_\omega \frac{p_j}{d_j}=\frac{\epsilon|}2$ from $(g_{i,j})^\omega=x_i$, $i=1,\ldots,k$. Therefore $z$ does not belong to any of the balls $B_C(\frac{\epsilon}{4},x_i)$. Hence $z$ does not belong to $B_C(\epsilon,\bar 1)$. But $\dist^\omega(z,\bar 1)=\lim_\omega \frac{\dist(z_j,1)}{d_j}\le \epsilon$, a contradiction.

Thus $\omega$-almost surely for every $j$, the ball $B_G(2p_j, 1)$ is contained in the union $$\bigcup_{i=1}^k g_{i,j}B_G(p_j,1).$$ But $B_G(2p_j,1)$ is equal to the product $B_G(p_j,1)B_G(p_j,1)$ in $G$. Hence $B_G(p_j, 1)$ is a $k$-approximable subgroup of $G$ $\omega$-almost surely for every $j$. Since $G$ is finitely generated and $\lim_\omega p_j=\infty$, every finite subset of $G$ is contained in one of the balls $B_G(p_j,1)$, so it remains to apply Theorem \ref{t:71}.
\endproof

{\bf Acknowledgement} This note was inspired by the paper \cite{Sis} where an attempt to answer the question of van den Dries and Wilkie was made. I am grateful to Terry Tao who suggested that the answer to that question should follow from \cite{H}.

\end{document}